\documentclass[12pt,twoside]{amsart}

\usepackage{amssymb}

\newtheorem{thm}{Theorem}

\newtheorem{cor}[thm]{Corollary}

\newtheorem{prop}[thm]{Proposition}

\theoremstyle{definition}
\newtheorem{defn}[thm]{Definition}

\newtheorem{say}[thm]{}

\newtheorem{rem}[thm]{Remark}          
\newtheorem{ack}{Acknowledgments}

\newtheorem{defn-thm}[thm]{Definition--Theorem}  %!!!!!!!!!!!!!!!!!!!!!!!!
\theoremstyle{remark}

\renewcommand{\c}[0]{{\mathbb C}}  

\renewcommand{\o}[0]{{\mathcal O}} 
\newcommand{\z}[0]{{\mathbb Z}}

\renewcommand{\a}[0]{{\mathbb A}}

\newcommand{\p}[0]{{\mathbb P}}

\newcommand{\q}[0]{{\mathbb Q}}

\newcommand{\qtq}[1]{\quad\mbox{#1}\quad}

\newcommand{\mult}[0]{\operatorname{mult}}

\newcommand{\diag}[0]{\operatorname{diag}}

\begin{document}
\bibliographystyle{amsplain}

\title[K\"ahler--Einstein metrics on
log del Pezzo surfaces]
{K\"ahler--Einstein metrics \\
on log del Pezzo surfaces \\
in weighted projective 3--spaces}
\author{Jennifer M.\ Johnson   and J\'anos Koll\'ar}

\maketitle

A {\it log del Pezzo} surface is a projective surface with
quotient singularities such that its anticanonical class is ample.
Such surfaces arise naturally in many different contexts, for instance
in connection with  affine surfaces \cite{miyanishi}, 
moduli of surfaces of general type \cite{alexeev94}, 3 and 4 dimensional
minimal model program \cite{alexeev93}.  They also provide a natural
testing ground for existence results of K\"ahler--Einstein metrics.
The presence of quotient singularities forces us to work with orbifold
metrics, but this is usually only a minor inconvenience.
Log del Pezzo surfaces with a K\"ahler--Einstein metric also
lead to Sasakian--Einstein 5--manifolds by \cite{bo-ga}.

In connection with \cite{dk}, 
the authors
 ran a computer program to find examples of
log del Pezzo surfaces in weighted projective spaces. The program
examined weights up to a few hundred and  produced 3 examples
of log del Pezzo surfaces where the methods of \cite[\S 6]{dk}
proved the existence
of a K\"ahler--Einstein metric.

The aim of this paper is twofold. First, we determine the complete
list of  anticanonically embedded quasi smooth 
log del Pezzo surfaces in weighted projective 3-spaces.
Second, we improve the methods of \cite[6.10]{dk}
to prove that many of these admit a K\"ahler--Einstein metric.
The same method also proves that some of these examples
do not have tigers (in the colorful terminology of \cite{keel-mc}). 

Higher dimensional versions of these results will be considered in
a subsequent paper.

\begin{defn}\label{wps.defn} For positive integers $a_i$
let $\p(a_0,a_1,a_2,a_3)$ denote the 
{\it  weighted
projective $3$-space} with weights $a_0\le a_1\le a_2\le a_3$.
(See \cite{dolg} or \cite{fletch}  for the basic definitons and results.)
We always assume that  any 3 of  the
  $a_i$ are relatively prime. 
We frequently write $\p$ 
to denote a weighted
projective $3$-space
if the weights are irrelevant or clear from the context.
We use $x_0,x_1,x_2,x_3$ to denote the corresponding weighted
projective coordinates. We let $(i,j,k,\ell)$ be an 
unspecified permutation of  $(0,1,2,3)$. 
$P_i\in \p(a_0,a_1,a_2,a_3)$ denotes the point
$(x_j=x_k=x_{\ell}=0)$.  The affine chart where
$x_i\neq 0$  can be written as
$$
\c^3(y_j,y_k,y_{\ell})/\z_{a_i}(a_j,a_k,a_{\ell}).
\eqno{(\ref{wps.defn}.1)}
$$
This shorthand denotes the quotient of $\c^3$ by the action
$$
(y_j,y_k,y_{\ell})\mapsto 
(\epsilon^{a_j}y_j,\epsilon^{a_k}y_k,\epsilon^{a_{\ell}}y_{\ell})
$$
where $\epsilon$ is a primitive $a_i$th root of unity.
The identification is given by $y_j^{a_i}=x_j^{a_i}/x_i^{a_j}$. 
(\ref{wps.defn}.1) are called the {\it orbifold charts} on
 $\p(a_0,a_1,a_2,a_3)$.

$\p(a_0,a_1,a_2,a_3)$ has  an index $a_i$ quotient singularity at $P_i$
and an index $(a_i,a_j)$ quotient singularity along the line
$(x_k=x_{\ell}=0)$.

For every $m\in \z$  there is a  rank 1 sheaf
$\o_{\p}(m)$ which is locally free only if $a_i|m$ for every $i$.
A basis of the space of  sections of $\o_{\p}(m)$
is given by all monomials in $x_0,x_1,x_2,x_3$ with weighted degree  $m$.
Thus $\o_{\p}(m)$ may have no sections for some $m>0$.
\end{defn}

\begin{say}[Anticanonically embedded quasi smooth surfaces]{\ }
\label{ant.emb.say}

Let $X\in |\o_{\p}(m)|$ be a surface of degree $m$. 
The adjunction formula
$$
K_X\cong \o_{\p}(K_{\p}+X)|_X \cong \o_{\p}(m-(a_0+a_1+ a_2+a_3))|_X
$$
holds iff  $X$ does not contain any of the singular lines.
If this condition holds then $X$ is a  (singular) del Pezzo surface
iff $m<a_0+a_1+ a_2+a_3$.  It is also well understood that 
from many points of view 
the most interesting
cases are when $m$ is as large as possible. Thus we consider the case
$X_d\in |\o_{\p}(d)|$ for $d=a_0+a_1+ a_2+a_3-1$.
We say that such an $X$ is {\it anticanonically embedded}.

Except for the classical cases 
$$
(a_0,a_1,a_2,a_3)= (1,1,1,1),\  (1,1,1,2)\qtq{or} (1,1,2,3),
$$
 $X$ is not smooth and  it passes through
some of the vertices $P_i$. Thus the best one can hope 
is that $X$ is smooth in the
orbifold sense, called {\it quasi smooth}. At the vertex $P_i$ this means that
the preimage of $X$ in the orbifold chart $\c^3(y_j,y_k,y_{\ell})$
is smooth. In terms of the equation of $X$ this is equivalent 
to saying that
$$
\mbox{For every $i$ there is a $j$ and 
 a monomial $x_i^{m_i}x_j$ of degree $d$.}
\eqno{(\ref{ant.emb.say}.1)}
$$
Here we allow $j=i$, corresponding to the case when
the general $X$ does not pass through $P_i$. 
The condition that $X$ does not contain any of the singular lines
is equivalent to
$$
\mbox{If $(a_i,a_j)>1$ then there is a monomial $x_i^{b_i}x_j^{b_j}$ 
of degree $d$.}
\eqno{(\ref{ant.emb.say}.2)}
$$
Finally, if every member of $|\o_{\p}(d)|$ contains a coordinate axis 
$(x_k=x_{\ell}=0)$ then the general member should be smooth along it, except
possibly at the vertices. That is
$$
\begin{array}{l}
\mbox{For every $i,j$, either  there is a monomial $x_i^{b_i}x_j^{b_j}$ 
of degree $d$,}\\
\mbox{or there are monomials  $x_i^{c_i}x_j^{c_j}x_k$
and $x_i^{d_i}x_j^{d_j}x_{\ell}$ of degree $d$.}
\end{array}
\eqno{(\ref{ant.emb.say}.3)}
$$
\end{say}

The computer search done in connection with
\cite{dk} looked at values of $a_i$ in a certain range
 to find the $a_i$ satisfying the constraints
(\ref{ant.emb.say}.1-3).  This approach starts with the $a_i$ and views
(\ref{ant.emb.say}.1-3) as linear equations in the unknowns $m_i, b_i,c_i,d_i$.
In order to find all solutions, we change the point of view.

\begin{say}[Description of the computer program]{\ }\label{descr.of.progr}

We
consider 
(\ref{ant.emb.say}.1) to be the main constraint,
 the $m_i$ as coefficients and the $a_i$ as
unknowns. The corresponding equations  can then be written
as a linear system 
$$
(M+J+U) (a_0\ a_1\ a_2\ a_3)^t
= (\ -1\ -1\ -1\ -1)^t
\eqno{(\ref{descr.of.progr}.1)}
$$
where $M=\diag(m_0,m_1,m_2,m_3)$ is a diagonal matrix,
$J$ is a matrix with all entries $-1$ and $U$ is
a matrix where each row has 3 entries $=0$ and one entry $=1$.
It is still not easy to decide when such a system has
positive integral solutions, but the main advantage is that
some of the $m_i$ can be bounded a priori.

Consider for instance $m_3$. The relevant equation is
$$
m_3a_3+a_j=a_0+a_1+a_2+a_3-1.
$$
Since $a_3$ is the biggest, we get right away that $1\leq m_3\leq 2$.
 Arguing inductively
with some case analysis we obtain that
\begin{enumerate}\setcounter{enumi}{1}
\item either  $2\leq m_2\leq 4 $ and $2\leq m_3\leq 10$,
\item or   the $a_i$ are in a series  $(1,a,b,b,)$ with $a|2b-1$. 
The latter satisfy (\ref{ant.emb.say}.2) only for $a=b=1$.
\end{enumerate}

Thus we have only finitely many possibilities for
the matrix $U$ and the numbers $m_1,m_2,m_3$.
Fixing these values, we obtain a linear system 
$$
(M+J+U) (a_0\ a_1\ a_2\ a_3)^t
= (\ -1\ -1\ -1\ -1)^t,
$$
where the only variable coefficient is the upper left corner
of $M$. Solving these formally we obtain that
$$
a_0=\frac{\gamma_0}{m_0\alpha+\beta}
$$
where $\alpha,\beta,\gamma_0$ depend only on $U$ and  $m_1,m_2,m_3$.
$a_0$ is supposed to be a positive integer, thus if
$\alpha\neq 0$ then there are only finitely many possibilities for $m_0$.
Once $m_0$ is also fixed, the whole system can be solved
and we check if the  $a_i$  are all positive integers.
We get 1362 cases.

If $\alpha=0$ but $\beta\neq 0$ then the general solution of
the system has the form
$$
a_0=\frac{\gamma_0}{\beta},
 a_i=\frac{m_0\delta_i+\gamma_i}{\beta}
\qtq{for} i=1,2,3.
$$
These generate the series of solutions, 405 of them.
Finally, with some luck, the case  $\alpha=\beta=\gamma_0=0$
never occurs, so we do not have to check further.

The resulting solutions need considerable cleaning up.
Many solutions 
$a_0,a_1,a_2,a_3$ occur multiply and we also have to check
the other conditions (\ref{ant.emb.say}.2-3).  
At the end we get the complete list, given in (\ref{main.thm}).

The computer programs are available at
\begin{verbatim}www.math.princeton.edu/~jmjohnso/LogDelPezzo
\end{verbatim}
\end{say}

These log del Pezzo surfaces are  quite interesting
in their own right. Namely, it turns out that
for many of them,
members of the linear systems $|-mK_X|$ can not be very 
singular at any point.  First we recall the notions log canonical etc.
(see, for instance,  \cite[2.3]{kmbook} for a detailed  introduction).

\begin{defn}\label{klt.etc.defn}
 Let $X$ be a surface and $D$ a $\q$-divisor on $X$.
Let $g:Y\to X$ be any proper birational morphism, $Y$ smooth.
 Then there is a unique
$\q$-divisor $D_Y=\sum e_iE_i$ on $Y$ such that 
$$
K_Y+D_Y\equiv g^*(K_X+D)\qtq{and}  g_*D_Y=D.
$$
We say that $(X,D)$ is
{\it canonical} (resp.\ {\it klt}, resp.\ {\it log canonical}) if
$e_i\geq 0$ (resp.\ $e_i>-1$, resp.\ $e_i\geq -1$) for every $g$
and for every $i$.
\end{defn}

\begin{defn} \cite{keel-mc}\label{tiger.defn}
 Let $X$ be a normal surface. A {\it tiger}
on $X$  is an effective  $\q$-divisor $D$ such that $D\equiv -K_X$ 
and $(X,D)$ is not klt.
As illustrated in \cite{keel-mc}, the tigers  carry
 important  information about birational transformations of 
log del Pezzo surfaces.
\end{defn}

\begin{rem} By a result of Shokurov (cf.\ \cite[22.2]{keel-mc}),
if the 
log del Pezzo surface $X$ has Picard number 1 and it has a tiger then
 $|\o_X(-mK_X)|\neq \emptyset$ for some $m=1,2,3,4,6$. 
The log del Pezzo surfaces in (\ref{main.thm}) mostly have bigger
 Picard number. It is quite interesting though that the two results
work for almost the same cases.
\end{rem}

We use the following sufficient condition to
obtain the existence of 
K\"ahler--Einstein metrics.

\begin{thm}\cite{nadel, dk}
\label{nadel.thm} Let $X$ be an $n$ dimensional Fano variety
(possibly with quotient singularities). Assume that
there is an $\epsilon>0$ such that 
$$
(X,\tfrac{n+\epsilon}{n+1}D) \qtq{is klt}
$$
for every effective $\q$-divisor $D\equiv -K_X$. Then
$X$ has a K\"ahler--Einstein metric.\qed
\end{thm}

The main result of this note is the following.

\begin{thm}\label{main.thm} There 
is  an  anticanonically embedded quasi smooth log del Pezzo surface
$X_d\subset \p(a_0,a_1,a_2,a_3)$ 
iff the $a_i$ and $d$ are among the following.
The table below also gives 
our results on the nonexistence of tigers (\ref{tiger.defn}) 
and on the existence of
K\"ahler--Einstein metrics. (Lower case 
$y$ means that the answer has been previously known.)
$$
\begin{array}{rcccccc}
a_0& a_1& a_2&a_3& d & \mbox{tiger} & \mbox{KE metric}\\
&&&&&&\\
\mbox{Series:}\quad 2&2k+1&2k+1& 4k+1& 8k+4& Y & Y\\
&&&&&&\\
\mbox{Sporadic:}\quad 1&1&1&1&3&y&y\\
1&1&1&2&4&y&y\\
1&1&2&3&6&y&y\\
1&2&3&5&10&y&?\\
1&3&5&7&15&y&?\\
1&3&5&8&16&y&?\\
2&3&5&9&18&?&?\\
3&3&5&5&15&N&Y\\
3&5&7&11&25&?&Y\\
3&5&7&14&28&?&Y\\
3&5&11&18&36&?&Y\\
5&14&17&21&56&N&Y\\
5&19&27&31&81&N&Y\\
5&19&27&50&100&N&Y\\
7&11&27&37&81&N&Y\\
7&11&27&44&88&N&Y\\
9&15&17&20&60&N&y\\
9&15&23&23&69&N&Y\\
11&29&39&49&127&N&Y\\
11&49&69&128&256&N&y\\
13&23&35&57&127&N&Y\\
13&35&81&128&256&N&y
\end{array}
$$
\end{thm}

\begin{rem}
The above results hold for every quasi smooth surface with the indicated
numerical data. 

Near the end of the list there are very few monomials 
of the given degree and in many cases there is only one such surface up
to isomorphism. In some other cases, for instance for the series,
 there are moduli.

 It is generally believed that the algebraic geometry
of any given log del Pezzo surface can be understood quite
well. 
There is every reason to believe that all of the remaining
 cases of (\ref{main.thm}) can be decided, though it may require a few pages 
of computation for each of them.
\end{rem}

\begin{say}[How to check if $(X,D)$ is klt or not?]{\ }
\label{how.to.check}

The definition (\ref{klt.etc.defn})
 requires understanding   all  resolutions of 
singularities. Instead, we use the following
multiplicity conditions  to check that a given divisor
is klt. These conditions are far from being necessary. 

Let $X$ be a surface with quotient singularities.
Let the singular points be $P_i\in X$ and we write these
locally analytically as 
$$
p_i:(\c^2, Q_i)\to (\c^2/G_i, P_i)\cong (X,P_i),
$$
 where $G_i\subset GL(2,\c)$
is a finite subgroup. We may assume that the origin is an
 isolated fixed point of every
nonidentity element of $G_i$ (cf. \cite{briesk}). 
Let $D$ be an effective
$\q$-divisor on $X$.
  Then  $(X,D)$ is klt if the following three
conditions  are satisfied.
\begin{enumerate}
\item (Non isolated non-klt points) $D$ does not contain an
irreducible component with coefficient  $\geq 1$.
\item (Canonical at smooth points) $\mult_PD\leq 1$ at every smooth point
$P\in X$. This follows from \cite[4.5]{kmbook}.
\item (Klt  at singular points) $\mult_{Q_i}D_i\leq 1$ for every $i$
where  $D_i:=p_i^*D$.
 This follows from \cite[5.20]{kmbook} and the previous case.
\end{enumerate}
In our applications we rely  on the following estimate.
\end{say}

\begin{prop}\label{WP.mult.est}
 Let $Z\subset \p(a_0,\dots, a_n)$ be a 
$d$-dimensional subvariety of a 
weighted projective space. 
Assume that $Z$ is not contained in the singular locus  and
that $a_0\leq \cdots\leq a_n$. Let
$Z_i\subset \a^n$ denote the preimage of $Z$ in the orbifold chart
$$
\a^n\to \a^n/\z_{a_i}\cong \p(a_0,\dots, a_n)\setminus (x_i=0).
$$
Then for every $i$ and every $p\in Z_i$,
$$
\mult_pZ_i\leq (a_n\cdots a_{n-d})(Z\cdot \o(1)^d).
$$
Moreover, if $Z\neq (x_0=\cdots=x_{n-d-1}=0)$ then we have a stronger
inequality
$$
\mult_pZ_i\leq (a_n\cdots a_{n-d+1}a_{n-d-1})(Z\cdot \o(1)^d).
$$
\end{prop}

Proof. Let $0\in C(Z)\subset \a^{n+1}$ denote the cone over $Z$
 with vertex $0$.
$Z_i$ can be identified with the hyperplane section $C(Z)\cap (x_i=1)$.
The multiplicty of a point is an upper semi continuous function on a variety,
thus it is sufficient to prove that
$$
\mult_0C(Z)\leq (a_n\cdots a_{n-d})(Z\cdot \o(1)^d).
$$
This is proved by induction on $\dim Z$. 

If $C(Z)$ is not contained
in the coordinate hyperplane $(x_i=0)$, then
write 
$$
Z\cap (x_i=0)=\sum_jm_jY_j\subset \p(a_0,\dots, a_n).
$$
Next we claim  that
\begin{eqnarray*}
\sum_jm_j(Y_j\cdot \o(1)^{d-1})& =& a_i(Z\cdot \o(1)^d)
\qtq{and}\\
 \sum_jm_j \mult_0 C(Y_j)& \geq  & \mult_0 C(Z).
\end{eqnarray*}
The first of these is the associativity of the intersection product,
and the second is a consequence of the usual estimate for the
intersection multiplicty (cf.\ \cite[12.4]{fulton}) applied to
$C(Z), (x_i=0)$ and $d-1$ other general hyperplanes through the origin.
(Note that in the first edition of \cite{fulton} there is a misprint in
(12.4). $\sum_{i=1}^r e_P(V_i)$ should be replaced by $\prod_{i=1}^r e_P(V_i)$.)
By the  inductive assumption   $\mult_0 C(Y_j)\leq (Y_j\cdot \o(1)^{d-1})$,
hence  $\mult_0 C(Z)\leq (Z\cdot \o(1)^d)$ as claimed.

In most cases, we can even choose  $i<n-d$. This is impossible only
if $Z\subset (x_0=\cdots=x_{n-d-1}=0)$, but then
equality holds.\qed

\begin{cor}\label{main.estimate.cor}
  Let $X_d\subset \p(a_0,a_1,a_2,a_3)$ be a quasismooth
surface of degree $d=a_0+a_1+a_2+a_3-1$. Then $X$ does not have a tiger
if $d\leq a_0a_1$.
If $(x_0=x_1=0)\not\subset X$ then $d\leq a_0a_2$
 is also sufficient.
\end{cor}

Proof.  Assume that $D\subset X_d$ is a tiger. We can view $D$
as a 1--cycle in $\p(a_0,a_1,a_2,a_3)$ whose degree is
$$
(D\cdot \o_X(1))=(\o_{\p}(d)\cdot \o_{\p}(1)\cdot \o_{\p}(1))
=\tfrac{d}{a_0a_1a_2a_3}.
$$
By (\ref{WP.mult.est}), this implies that
the multiplicity of $D_i$ (as in (\ref{how.to.check}.3)) 
is bounded from above by $\tfrac{d}{a_0a_1}$ at any point.
 Thus $(X,D)$ is klt if $d\leq a_0a_1$.

If $(x_0=x_1=0)\not\subset X$ then we can weaken this to $d\leq a_0a_2$,
again by (\ref{WP.mult.est}).\qed

\medskip

Using (\ref{nadel.thm}) and  a similar argument we  obtain the following.

\begin{cor}\label{KE.estimate.cor}
  Let $X_d\subset \p(a_0,a_1,a_2,a_3)$ be a quasismooth
surface of degree $d=a_0+a_1+a_2+a_3-1$. Then $X$ admits a
K\"ahler--Einstein metric if
 $d< \tfrac32 a_0a_1$.
If $(x_0=x_1=0)\not\subset X$ then $d< \tfrac32  a_0a_2$
 is also sufficient.\qed
\end{cor}

\begin{say}[Proof of (\ref{main.thm})]{\ }

The nonexistence of tigers and the existence of a K\"ahler--Einstein metric
in the 
sporadic examples follows from
(\ref{main.estimate.cor}) and 
(\ref{KE.estimate.cor}).
There are 5 cases when we need to use that $X$ does not contain the
line $(x_0=x_1=0)$. This is equivalent to claiming that
the equation of $X$ contains a monomial involving $x_2,x_3$ only.
In all 5 cases this is already forced by the condition (\ref{ant.emb.say}.1).

Assume next that  $X$ is one of the series $(2,2k+1,2k+1,4k+1)$. 
Its equation is a linear combination of terms
$$
x_0^{4k+2}, x_3^2x_0, x_3(x_1+x_2)x_0^{k+1},g_4(x_1,x_2), 
g_2(x_1,x_2)x_0^{2k+1}.
$$
Moreover, the conditions (\ref{ant.emb.say}.1--3) imply that
the first 2 appear with nonzero coefficient and $g_4$ does
 not have multiple roots.

$(x_0=0)$ intersects $X$ in  a curve $C$ whith equation
$$
\begin{array}{l}
(q_{8k+4}(x_1,x_2)=0)\subset \p^2(2k+1,2k+1,4k+1),
\qtq{ isomorphic to}\\
 (q_{4}(x_1,x_2)=0)\subset \p^2(1,1,4k+1).
\end{array}
$$
Thus $C$ has 4 irreducible components $C_1+C_2+C_3+C_4$
meeting at $P_3=(0:0:0:1)$.
This shows that $\frac12 C$ is not klt at $P_3$ and 
$\frac12 C$ is a tiger on $X$.

Next we prove that $(X,D)$ is log canonical for every 
effective $\q$-divisor $D\equiv -K_X$. This is stronger than needed
in order to apply (\ref{nadel.thm}).

Consider the linear system $\o_{\p}(2(2k+1))$. This is the pull back
of $\o(2(2k+1))$ from the weighted projective plane $\p(2,2k+1,2k+1)$.
The latter  is isomorphic to $\p(2,1,1)$ which is the quadric cone in 
ordinary $\p^3$ and the linear system is the hyperplane sections,
thus very ample. Hence for every smooth point $P\in X$
there is a divisor $F\in |\o_{X}(2(2k+1))|$ passing through $P$
and not containing any of the irreducible components of $D$.
So
$$
\mult_PD\leq (D\cdot F)=\frac{2(2k+1)(8k+4)}{2(2k+1)^2(4k+1)}=
\frac{4}{4k+1}<1.
$$
We are left to deal with the singular points of $X$. These are
at $P_3=(0:0:0:1)$ and at $P_a=(0:a:1:0)$ where $a$ is a root
of $g_4(x_1,1)$.

$P_3$ is the most interesting. 
Let $p_3: (S\cong \c^2, Q_3)\to (X,P_3)$ be a local orbifold chart.
Intersecting  $p_3^*D$ with a general member of 
the linear system $|x_0^{2k+1}, x_1^2|$  we obtain that
$$
\mult_{Q_3} p_3^*D\leq \tfrac{4k+1}{2} (D\cdot \o(2(2k+1))=2.
$$
This is too big to apply (\ref{how.to.check}.3). Let
$\pi:S'\to S$ be the blow up of the origin with exceptional divisor
$E$. Then
$$
K_{S'}+\alpha E+\pi^{-1}_*(p_3^*D)\equiv \pi^*(K_S+p_3^*D),
$$
and $\alpha\leq 1$. Using Shokurov's inversion of adjunction
 (see, for instance \cite[5.50]{kmbook}) 
$(X,D)$ is log canonical at
$P_3$ if $\pi^{-1}_*(p_3^*D)|_E$ is a sum of points, all with coefficient
$\leq 1$.  In order to estimate these coefficients, we
write 
$D=D'+\sum a_iC_i$ where $D'$ does not contain any of the $C_i$.

We first compute that
$$
(C_i\cdot C_j)=\tfrac1{4k+1}\qtq{if $i\neq j$, and}
(C_i\cdot \o(1))=\tfrac{1}{(2k+1)(4k+1)}.
$$
From this we obtain that 
$$
(C_i\cdot C_i)=(C_i\cdot \o(1))-\sum_{j\neq i}(C_i\cdot C_j)=
\tfrac{-(6k+1)}{(2k+1)(4k+1)}.
$$
Thus
$$
\tfrac{1}{(2k+1)(4k+1)}=(C_i\cdot D)=
a_i(C_i\cdot C_i)+(\textstyle{\sum_{j\neq i}a_j})(C_i\cdot C_{i+1})+
(C_i\cdot D')
$$
Multiplying by $(2k+1)(4k+1)$ and
using that $\sum a_i\leq  \mult_{Q_3} p_3^*D\leq 2$
and $(C_i\cdot D')\leq (D\cdot D)$ this becomes
$$
1\leq -(6k+1)a_i+(2-a_i)(2k+1)+4
\qtq{which gives}  a_i\leq \tfrac12+\tfrac{2}{4k+1}.
$$
Furthermore,
$$
\mult_{Q_3} p_3^*D'\leq \tfrac14(p_3^*D'\cdot \sum_ip_3^*C_i)
\leq \tfrac{4k+1}{4}(D\cdot \o(2))=\tfrac1{2k+1}.
$$
Thus we see that 
$$
\pi^{-1}_*(p_3^*D)|_E=\sum a_i\pi^{-1}_*(p_3^*C_i)|_E+
\pi^{-1}_*(p_3^*D')|_E
$$
is a sum of 4 distinct points with coefficient $\leq \tfrac12+\tfrac{2}{4k+1}$
and another sum of points where the sum of the coefficients is
 $\leq \tfrac1{2k+1}$. Since $\tfrac12+\tfrac{2}{4k+1}+\tfrac1{2k+1}< 1$,
we see that $(X,D)$ is log canonical at $P_3$.

The points $P_a$ are easier. Only one of the $C_i$
passes through each of them, and the multiplicity of
the pull back of $D'$ is bounded by 
$\tfrac{2k+1}{4}(D\cdot \o(2))=\tfrac1{4k+1}$. This shows right away
that $(X,D)$ is klt at these points.\qed
\end{say}

\begin{ack}  We   thank  Ch.\ Boyer, J.-P.\ Demailly  and S.\ Keel
for  helpful
comments and references.
Partial financial support was provided by  the NSF under grant number 
DMS-9970855. 
\end{ack}

\vskip1cm

\noindent Princeton University, Princeton NJ 08544-1000

\begin{verbatim}jmjohnso@math.princeton.edu\end{verbatim}

\begin{verbatim}kollar@math.princeton.edu\end{verbatim}
\end{document}